\documentclass[preprint]{elsarticle}
\makeatletter
\def\ps@pprintTitle{%
 \let\@oddhead\@empty
 \let\@evenhead\@empty
 \def\@oddfoot{\centerline{\thepage}}%
 \let\@evenfoot\@oddfoot}
\makeatother

\usepackage{amsmath}
\usepackage{amssymb}

\newtheorem{thr}{Theorem}
\newtheorem{lmm}{Lemma}
\newtheorem{prp}{Proposition}
\newtheorem{clm}{Claim}
\newtheorem{obs}{Observation}
\newtheorem{cnj}{Conjecture}
\newdefinition{dfn}{Definition}
\newdefinition{case}{Case}
\newproof{prf}{Proof}

\begin{document}

\begin{frontmatter}

\title{A New Proof for the Characterization of Linear Betweenness Structures}

\author[arim]{P\'eter G. N. Szab\'o\corref{cr1}}
\cortext[cr1]{Corresponding author.}
%\tnotetext[fn1]{This work was supported by the National Research, Development and Innovation Office -- NKFIH, No. 108947.}

\ead[arim]{szape@cs.bme.hu}
%\date{\today} % nem csina'l semmit

\address[arim]{Department of Computer Science and Information Theory, Budapest University of Technology and Economics, M\H{u}egyetem rkp. 3., H-1111 Budapest, Hungary}

\begin{abstract}
In their paper published in 1997, Richmond and Richmond classified metric spaces in which all triangles are degenerate. That result was later reproved by Dovgoshei and Dordovskii in the finite case and it was generalized to finite pseudometric betweennesses by Beaudou et al. In this paper, we give a new, independent proof to the finite case of the original theorem which we reformulate in terms of linearity of betweenness structures.
\end{abstract}


\end{frontmatter}

\section{Introduction}
One can say that metric space is one of the most successful concepts of mathematics, with various applications in several fields including, among others, computer science, quantitative geometry, topology, molecular chemistry and phylogenetics.
Although finite metric spaces are trivial objects from a topological point of view, they have surprisingly complex and intriguing combinatorial properties which were investigated from different angles over the last fifty years. The concept of metric betweenness appears to play a central role in the related literature which span from combinatorial geometry to metric graph theory.

The largest field of related research, metric graph theory, studies different classes of graphs in terms of the induced graph metric and the geodesic betweenness.
One of the earliest results are Sholander's axiomatic characterization of trees, lattices and partially ordered sets in terms of segments, medians and betweenness \cite{sholander1952trees}.
The (geodesic) interval function of a connected graph was first extensively studied by Mulder \cite{mulder1980interval}, who introduced the five classical axioms of the interval function.
In a series of papers, Nebesk\'{y} et al. described the interval function of a connected graph in terms of first order transit axioms, each time improving the proof \cite{nebesky1994characterization, nebesky1994characterization-set, nebesky1998characterizing, nebesky2001characterization, nebesky2001interval, mulder2009axiomatic}.
Besides the interval function, two other path transit functions were extensively studied on connected graphs: the induced path function and the all-paths function. Mulder introduced a general notion of transit function \cite{mulder2008transit} to unify the three concepts and presented a list of prototype problems that were only studied for specific transit functions but are unsolved in the general case.
For a thorough survey on geodesic and induced path betweenness, see \cite{changat2019betweenness}.

There are two other important lines of research related to the study of metric betweenness.
Based on the pioneering work of Isbell \cite{isbell1964six} and Buneman \cite{buneman1974note}, Dress et al. studied both algorithmic and combinatorial aspects of phylogenetic trees and the split decomposition of finite metric spaces \cite{bandelt1992canonical, dress1987parsimonious}. These results have important applications in evolutionary biology.
%Another purely combinatorial approach was presented by Mascioni \cite{mascioni2004equilateral}, who proved novel Ramsey-type results for finite metric spaces.
A more recent open problem is the generalization of the de Bruijn--Erd\H{o}s theorem to finite metric spaces, originally conjectured by Chen and Chv\'atal in \cite{chen2008problems}.
We note that the particular definition of line used there is essentially different from the one we introduce in this paper.
The conjecture is still open today, however, it has been already proved in a number of important cases: for some subspaces of the Euclidean plane with $L_1$ and $L_\infty$ metric \cite{kantor2013debruijn}, for finite $1$-$2$ metric spaces \cite{chvatal2014debruijn}, for chordal and distance-hereditary graphs \cite{aboulker2018new}, for bisplit graphs \cite{beaudou2019bisplit} and for $(q, q - 4)$-graphs \cite{schrader2019debruijn}. Further, polynomial lower bounds have been proved in the general case of pseudometric, metric and graphic betweennesses \cite{aboulker2016lines}.

In \cite{richmond1997metric}, Richmond and Richmond obtained a nice characterization of metric spaces which does not contain degenerate triangles, i.e. triangles where the sum of two sides is equal to the third side. This result was later reproved by Dovgoshei and Dordovskii in the finite case \cite{dovgoshei2009betweenness} and generalized to finite pseudometric betweennesses by Beaudou et al. \cite{beaudou2013lines}.

In this paper, we give a new, independent proof to the above theorem of Richmond and Richmond in the finite case but we discuss it from the perspective of linearity of betweenness structures, which is equivalent to the property of having no degenerate triangles.
First, we introduce the framework and system of notations that will be used throughout the paper.
A \emph{metric space} is a pair $M = (X, d)$ where $X$ is a nonempty set and $d$ is a \emph{metric} on $X$, i.e. an $X\times X\rightarrow\mathbb{R}$ function which satisfies the following conditions for all $x, y, z\in X$:
\begin{enumerate}
\item $d(x, y) = 0\Leftrightarrow x = y$ (\emph{identity of indiscernibles})\label{Ems1};
\item $d(x, y) = d(y, x)$ (\emph{symmetry})\label{Ems2};
\item $d(x, z)\leq d(x, y) + d(y, z)$ (\emph{triangle inequality})\label{Ems3}.
\end{enumerate}
The non-negativity of metric follows from the definition.
We will refer to the ground set and the metric of the metric space $M$ by $X(M)$ and $d_M$, respectively.

If the triangle inequality holds with equality for three points $x, y, z\in X$, i.e. $d(x, z) = d(x, y) + d(y, z)$, we write $(x\ y\ z)_M$ (or simply $(x\ y\ z)$ if $M$ is clear from the context), and we say that $y$ is \emph{between} $x$ and $z$ in $M$. We call this ternary relation the betweenness relation of $M$. Further, if $(x\ y\ z)$ holds, we say that $x$, $y$ and $z$ are \emph{collinear}.
In the rest of the paper, every metric space will be assumed to be \emph{finite} ($|X(M)| <\infty$) if not stated otherwise.

The relation of betweenness of a metric space has the following elementary properties.
For all $x, y, z\in X$,
\begin{enumerate}
\item $(x\ x\ z)$\label{Ecoll2};
\item $(x\ y\ z)\Rightarrow (z\ y\ x)$\label{Ecoll1};
\item $(x\ y\ z)\wedge (y\ x\ z)\Rightarrow x = y$\label{Ecoll3}.
\end{enumerate}
The \emph{trichotomy} of betweenness follows straight from these properties:
for any three distinct points $x, y, z\in X$, at most one of the relations $(x\ y\ z)$, $(y\ z\ x)$, $(z\ x\ y)$ hold.

Different metrics on ground set $X$ may define the same betweenness relation. Since we are interested in the combinatorial properties of the betweenness relation, we do not need to know the exact values of the underlying metric. Therefore, we base our definitions and theorems on the abstraction level of so-called betweenness structures as described below.

A \emph{betweenness structure} is a pair $\mathcal{B} = (X,\beta)$ where $X$ is a nonempty finite set and $\beta\subseteq X^3$ is a ternary relation, called the \emph{betweenness relation} of $\mathcal{B}$.
The fact $(x, y, z)\in\beta$ will be denoted by $(x\ y\ z)_\mathcal{B}$ or simply by $(x\ y\ z)$ if $\mathcal{B}$ is clear from the context, and we say that $x$, $y$ and $z$ are \emph{collinear} and that $y$ is \emph{between} $x$ and $z$.

The \emph{substructure} of $\mathcal{B}$ induced by a nonempty subset $Y\subseteq X$ is the betweenness structure $\mathcal{B}\vert_Y = (Y,\beta\cap Y^3)$. The substructure $\mathcal{B}\vert_{X\backslash\{x\}}$ will also be denoted by $\mathcal{B} - x$.
Two betweenness structures $\mathcal{B} = (X,\beta)$ and $\mathcal{C} = (Y,\gamma)$ are \emph{isomorphic} (in notation $\mathcal{B}\simeq\mathcal{C}$) if there exists a bijection $\varphi: X\rightarrow Y$ such that for all $x, y, z\in X$, $(x\ y\ z)_\mathcal{B}\Leftrightarrow(\varphi(x)\ \varphi(y)\ \varphi(z))_\mathcal{C}$.

There is a natural way to associate a betweenness structure with a metric space: the \emph{betweenness structure induced by a metric space} $M = (X, d)$ is $\mathcal{B}(M) = (X,\beta_M)$ where $\beta_M$ is the betweenness relation of $M$, as defined above.
To simplify notations, we will write $(x\ y\ z)_M$ for $(x\ y\ z)_{\mathcal{B}(M)}$.

A betweenness structure is \emph{metric} if it is induced by a metric space.
We note that the same elementary properties hold for the betweenness relation of a metric betweenness structure that hold for the betweenness relation of a metric space, including trichotomy.
Further, substructures of a metric betweenness structure are metric as well.
In the rest of the paper, every betweenness structure will be assumed to be metric if not stated otherwise.

\begin{figure}[ht]
\centering
\includegraphics{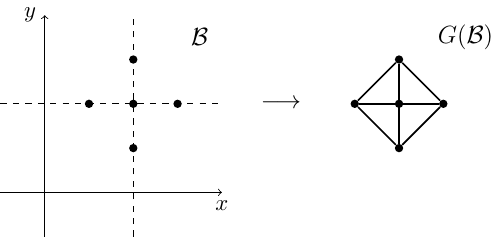}
\caption{The underlying graph of a betweenness structure defined by five points on the Euclidean plane}\label{Fadj}
\end{figure}

By graph we always mean a simple graph.
The \emph{underlying graph} (or adjacency graph) of a betweenness structure $\mathcal{B} = (X,\beta)$ is the graph $G(\mathcal{B}) = (X, E(\mathcal{B}))$ where the edges are such pairs of distinct points for which no third point lies between them (see Figure \ref{Fadj}). More formally,
$$E(\mathcal{B}) =\left\{\{x,z\}\in\binom{X}{2}:\nexists\,y\in X\backslash\{x, z\},\,(x\ y\ z)_\mathcal{B}\right\}.$$
These edges are sometimes called primitive pairs in the related literature.
The underlying graph is our most important connection to graph theory. Not only it is a helpful tool in examining the underlying betweenness structure but it is the ``minimal'' graph that can induce the underlying metric space with an appropriate edge weighting.

Let $G = (V, E)$ be a connected graph. The \textit{metric space induced by $G$} is the metric space $M(G) = (V, d_G)$ where $d_G$ is the usual \emph{graph metric} of $G$, i.e. $d_G(u, v)$ is the length of the shortest path between $u$ and $v$ in $G$.
The \emph{betweenness structure induced by $G$} is the betweenness structure induced by $M(G)$, also denoted by $\mathcal{B}(G)$.
It is easy to see that $(x\ y\ z)_{\mathcal{B}(G)}$ holds if and only if $y$ is on a shortest path connecting $x$ and $z$ in $G$.
A betweenness structure (or a metric space) is \emph{graphic} if it is induced by a connected graph. We note that any connected graph $G$ satisfies $G(\mathcal{B}(G)) = G$, however, $\mathcal{B}(G(\mathcal{B})) = G$ is true if and only if $\mathcal{B}$ is graphic.

\section{Linear Betweenness Structures}

In this section, we state and discuss the central result of this paper, the characterization of linear betweenness structures (Theorem \ref{T1}). We also compare our definition of line to the one used by Chen and Chv\'atal in \cite{chen2008problems}.

As usual, $P_n$ and $C_n$ denotes the path and the cycle of length $n$, respectively. Additionally, we assume for convenience that the set of vertices of $P_n$ and $C_n$ are the integers from $1$ to $n$ and the edges are the pairs of consecutive integers (where $n$ and $1$ are considered consecutive in case of $C_n$).
Let $\mathcal{P}_n$ and $\mathcal{C}_n$ denote the graphic betweenness structures induced by $P_n$ and $C_n$, respectively.
A betweenness structure is \emph{ordered} if it is induced by a path, or in other words, it is isomorphic to $\mathcal{P}_n$. Such an isomorphism is called an \emph{ordering}. We denote the ordered betweenness structure induced by the path $P = x_1x_2\ldots x_n$ by $[x_1, x_2,\ldots, x_n]$.

Let $\mathcal{B} = (X,\beta)$ be a betweenness structure.
We say that a set $Y\subseteq X$ is \textit{collinear} in $\mathcal{B}$ if any three points of $Y$ are collinear.
A \textit{line} of $\mathcal{B}$ is a maximal collinear set of $\mathcal{B}$.

\begin{dfn}
A betweenness structure $\mathcal{B} = (X,\beta)$ is \textit{linear} if any three points of $\mathcal{B}$ are collinear, i.e. $X$ is itself a line of $\mathcal{B}$.
\end{dfn}

Observe that all ordered betweenness structures are linear and substructures of a linear betweenness structure are linear as well.
Since every line induces a linear substructure and the ground set of a linear betweenness structure is itself a line, the two definitions describe the same concept from different points of view.

It is a natural question to ask, what the linear betweenness structures are up to isomorphism, and further, what is the difference between linearity and orderedness. The answer was first given by Richmond and Richmond in \cite{richmond1997metric}, albeit in a slightly different form. We reformulate that result to match our definitions.

\begin{thr}[Richmond, Richmond \cite{richmond1997metric}]\label{T1}
A betweenness structure $\mathcal{B}$ is linear if and only if $\mathcal{B}\simeq\mathcal{P}_n$ (for $n\geq 1$) or $\mathcal{B}\simeq\mathcal{C}_4$.
\end{thr}

What Theorem \ref{T1} tells us, maybe a little surprisingly, is that there is exactly one unordered linear betweenness structure up to isomorphism, which also means that linearity is not equivalent to the much simpler property of orderedness but only by a hair's breadth.

Based on Theorem \ref{T1}, we can divide lines into two subcategories. The ones isomorphic to $\mathcal{P}_n$ (for some $n\geq 1$) we call \emph{ordered line}s, while the ones isomorphic to $\mathcal{C}_4$ we call \emph{cyclic line}s. We say that a betweenness structure is \emph{regular} if it is $\mathcal{C}_4$-free i.e. it does not contain a cyclic line as a substructure.
The usefulness of this distinction is justified by several problems we encountered that are much easier to solve for regular betweenness structures than for irregular ones.

Before we present the new proof to Theorem \ref{T1}, we want to point out some of the similarities and differences between our definition of line and the one used in \cite{chen2008problems} that will be respectively called \emph{tight} line and \emph{spanned} line in order to avoid ambiguity.
A \emph{spanned line} in a metric space $M$ (or in a betweenness structure $\mathcal{B}$) on $X$ is a subset of $X$ of the form $\overline{xy} =\{z\in X: x, y, z\text{ are collinear in } M\}$ where $x$ and $y$ are distinct points in $X$.

The two definitions of line have different advantages and disadvantages as they preserve different key properties of the Euclidean line.
On one hand, there is a natural way to associate a spanned line with a pair of distinct points and thus there can be at most $\binom{|X|}{2}$ of them. This is obviously not true for tight lines. On the other hand, however, tight lines --unlike spanned lines-- cannot contain each other, hence, form a Sperner system on $X$.

A spanned line is called \emph{universal} if it is equal to $X$. In \cite{debruijn1948combinatorial}, de Bruijn and Erd\H{o}s showed that $n$ distinct points on the Euclidean plane determine at least $n$ distinct lines.
Interest towards this result were renewed in 2008, when Chen and Chv\'atal conjectured that it may generalize to finite metric spaces as follows.

\begin{cnj}[Chen, Chv\'atal \cite{chen2008problems}]\label{Cchvatal}
If a betweenness structure $\mathcal{B}$ does not contain a universal spanned line, then there are at least $n$ distinct spanned lines in $\mathcal{B}$.
\end{cnj}

We can say that a betweenness structure is linear in the ``spanned'' sense if it contains a universal line. The two definitions of linearity are quite similar in the respect that both requires $X$ to be covered by a line, although, our definition of linearity is more restrictive.
However, there are several essential differences. Most importantly, Theorem \ref{T1} shows that tight lines are better for generalizing orderedness, as spanned lines can be very far from being ordered. The de Bruijn--Erd\H{o}s theorem, however, would not generalize so nicely with tight lines, as the following proposition shows.

\begin{prp}\label{P3lines}
If $\mathcal{B}$ is a non-linear betweenness structure of order $n$, then there are at least $3$ tight lines in $\mathcal{B}$ and this bound is best possible.
\end{prp}
\begin{prf}
On one hand, if $\mathcal{B} = (X,\beta)$ is non-linear, then there exist three distinct points in $X$ that are not collinear. These points obviously determine $3$ distinct tight lines in $\mathcal{B}$. On the other hand, every betweenness structure induced by a tree with exactly $3$ leaves has exactly $3$ tight lines. $\square$
\end{prf}

\section{A New Proof of Theorem \ref{T1}}
First, we list some general helper statements that we will use in the proof later on.
Observation \ref{Ofour} is a well-known property of metric betweenness, while Observation \ref{Oconn} and \ref{Osubmetr} encapsulate simple properties of the underlying graph.
Lastly, we reformulate a particularly interesting remark of Dress (Remark 3 in \cite{dress2007category}) as Proposition \ref{Pdress}.
We will only use it in the special case when $T$ is a path.

\begin{obs}[Four Relations]\label{Ofour}
Let $\mathcal{B}$ be a betweenness structure on $X =\{x_1, x_2, x_3, x_4\}$ such that $(x_1\ x_3\ x_4)$ and $(x_1\ x_2\ x_3)$ hold. Then $(x_1\ x_2\ x_4)$ and $(x_2\ x_3\ x_4)$ hold as well.
\end{obs}

\begin{obs}\label{Oconn}
The underlying graph of a betweenness structure is connected.
\end{obs}

\begin{obs}\label{Osubmetr}
Let $\mathcal{B}$ be a betweenness structure and let $Y$ be a nonempty set of points in $\mathcal{B}$. Then $G(\mathcal{B})[Y]\leq G({\mathcal{B}\vert_Y})$.
\end{obs}

%%% prove for almost-metric bs's - postpone for later
\begin{prp}[Dress \cite{dress2007category}]\label{Pdress}
Let $\mathcal{B}$ be a betweenness structure such that\\ $G(\mathcal{B}) = T$ is a tree. Then $\mathcal{B}$ is induced by $T$.
\end{prp}

Now, let $\mathcal{B} = (X,\beta)$ be a betweenness structure as in Theorem \ref{T1}, let $n = |X|$ and $G = G(\mathcal{B})$.
First, we note that the ``if'' part of the theorem is obvious because both $\mathcal{C}_4$ and $\mathcal{P}_n$ are clearly linear. We also note that the ``only if'' part trivially holds for $n\leq 3$. 
Hence, it remains to show that
\begin{itemize}
\item if $n = 4$, then either $\mathcal{B}\simeq\mathcal{P}_4$ or $\mathcal{B}\simeq\mathcal{C}_4$;
\item if $n\geq 5$, then $\mathcal{B}\simeq\mathcal{P}_n$.
\end{itemize}

First, we prove two main lemmas and then proceed with a smallest counterexample argument.
The first lemma is the key tool to prove that $\mathcal{B}$ is induced by a path or a cycle of length $4$.

\begin{lmm}\label{Lpathimpl}
If $G\simeq P_n$ ($n\geq 1$) or $G\simeq C_4$, then $\mathcal{B} =\mathcal{B}(G)$.
\end{lmm}
\begin{prf}
Case $G\simeq P_n$ follows straight from Proposition \ref{Pdress} with $T = G$.
Suppose now that $G\simeq C_4$ and let $x_1, x_2, x_3, x_4$ denote the vertices of $G$ (in consecutive order). The next observation follows from the linearity of $\mathcal{B}$.

\begin{obs}\label{Oedges}
Let $x, y, z\in X$ be distinct points such that both $x$ and $z$ are adjacent to $y$ in $G$. Then $(x\ y\ z)$ holds.
\end{obs}

Now, $(x_1\ x_2\ x_3), (x_2\ x_3\ x_4), (x_3\ x_4\ x_1)$ and $(x_4\ x_1\ x_2)$ follows from Observation \ref{Oedges}, which means exactly that $\mathcal{B} =\mathcal{B}(G)$.
$\square$
\end{prf}

From Lemma \ref{Lpathimpl} we obtain immediately that if $G\simeq P_n$ or $G\simeq C_4$, then $\mathcal{B}\simeq\mathcal{P}_n$ or $\mathcal{B}\simeq\mathcal{C}_4$, respectively. Our goal in the rest of the proof is thus to prove that $G\simeq P_n$ or $G\simeq C_4$.

\begin{lmm}\label{Lgrlin}
The graph $G$ is either a path or a cycle.
\end{lmm}
\begin{prf}
Observation \ref{Oconn} shows that $G$ is connected, therefore, it is enough to show that $\Delta(G)\leq 2$.
Assume to the contrary that there exists a point $x\in X$ such that $d_G(x)\geq 3$.
Let $u, v$ and $w$ be three distinct neighbors of $x$. Because of linearity, we can assume, for example, that $(u\ v\ w)$ holds.
Now, $(u\ x\ v)$ and $(v\ x\ w)$ hold by Observation \ref{Oedges}. 
However, by applying Observation \ref{Ofour} to $(u\ v\ w)$ and $(u\ x\ v)$, we obtain $(x\ v\ w)$ in contradiction with trichotomy. $\square$
\end{prf}

Now, by virtue of Lemma \ref{Lgrlin}, we only need to prove that $G\not\simeq C_n$ if $n\geq 5$.
Suppose to the contrary that this statement is false. We can assume that $\mathcal{B}$ is a smallest counterexample. Let $\varphi$ be an isomorphism between $G$ and $C_n$, $x_i =\varphi^{-1}(i)$ and $G_i = G(\mathcal{B} - x_i)$ for all $i\in [n]$.

\begin{figure}[ht]
\centering
\includegraphics{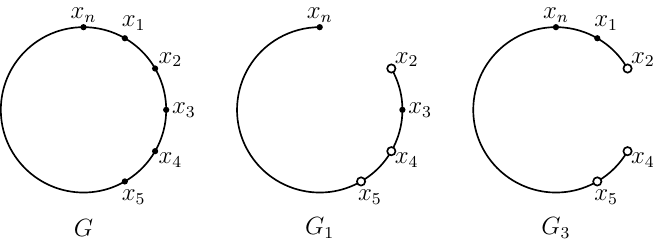}
\caption{Proof of Claim \ref{Ctwopath}, case $n\geq 5$, $G_1\simeq G_3\simeq P_{n - 1}$}\label{Ftwopath1}
\end{figure}

\begin{figure}[ht]
\centering
\includegraphics{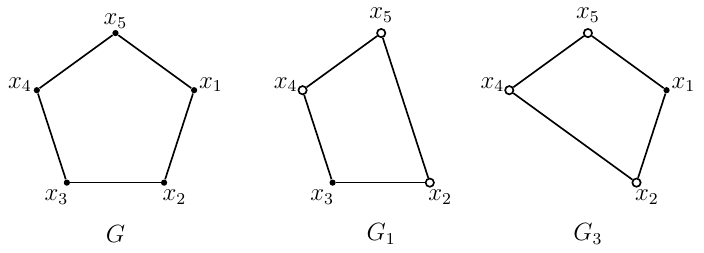}
\caption{Proof of Claim \ref{Ctwopath}, case $n = 5$, $G_1\simeq G_3\simeq C_4$}\label{Ftwopath2}
\end{figure}

\begin{clm}\label{Ctwopath}
Let $x_i$ and $x_j$ be two vertices of $G$ at distance $2$. Then 
\begin{enumerate}
\item $G_i\simeq P_{n - 1}\Rightarrow G_j\not\simeq P_{n - 1}$;\label{Etwopath1}
\item $G_i\simeq C_4\Rightarrow G_j\not\simeq C_4$.\label{Etwopath2}
\end{enumerate}
\end{clm}
\begin{prf}
Without loss of generality, we can assume that $i = 1$ and $j = 3$.
First, suppose to the contrary that $G_1\simeq G_3\simeq P_{n - 1}$.
It is obvious from Observation \ref{Osubmetr} that $G_1 = G - x_1$. Now, we can apply Proposition \ref{Pdress} with $T = G_1$ to obtain that $\mathcal{B} - x_1 =\mathcal{B}(G_1) =\mathcal{B}(G - x_1) = [x_2, x_3,\ldots, x_n]$ and consequently $(x_2\ x_4\ x_5)$ holds (see Figure \ref{Ftwopath1}).
Similarly, $\mathcal{B} - x_3 =\mathcal{B}(G_3) =\mathcal{B}(G - x_3) = [x_4, x_5,\ldots, x_n, x_1, x_2]$, therefore, $(x_4\ x_5\ x_2)$ holds, which contradicts the previous assertion.

For Part \ref{Etwopath2}, suppose that $G_1\simeq G_3\simeq C_4$.
Again, it is obvious from Observation \ref{Osubmetr} that the edges of $G_1$ are $\{x_2,x_3\}$, $\{x_3,x_4\}$, $\{x_4,x_5\}$ and $\{x_5,x_2\}$ (see Figure \ref{Ftwopath2}). Since $\mathcal{B} - x_1$ is a linear betweenness structure as well, we obtain that $\mathcal{B} - x_1 =\mathcal{B}(G_1)$ by Lemma \ref{Lpathimpl}, which further implies $(x_2\ x_5\ x_4)$.
Similarly, if $G_3\simeq C_4$, then $(x_2\ x_4\ x_5)$ would hold in contradiction with trichotomy. $\square$
\end{prf}

Now, the following two cases complete the proof.

\begin{case}[$n = 5$]
Since $\mathcal{B} - x_i$ is a linear betweenness structure, $G_i$ is isomorphic to either $P_4$ or $C_4$ by Lemma \ref{Lgrlin}. Color the vertex $i$ of $C_5$ \emph{red} if $G_i\simeq P_4$, \emph{blue} otherwise.
Because of Claim \ref{Ctwopath}, non-adjacent vertices are colored differently, hence, we obtain a proper $2$-vertex-coloring of $\overline{C}_5\simeq C_5$, which is obviously a contradiction.
\end{case}

\begin{case}[$n > 5$]
Since for all $i\in [n]$, $\mathcal{B} - x_i$ is a linear betweenness structure of at least five points, the minimality of $\mathcal{B}$ and Lemma \ref{Lgrlin} implies that $G_i\simeq P_{n - 1}$ in contradiction with Claim \ref{Ctwopath}. $\square$
\end{case}

\section{Conclusion}

In \cite{richmond1997metric}, Richmond and Richmond showed that a metric space without degenerate triangles can be isometrically embedded into the real line with only one exception. In this paper, we reformulated this result using the notion of linearity of betweenness structures and presented a new, independent proof to it in the finite case.

As for future research, we plan to study the extremal cases of Proposition \ref{P3lines} in a forthcoming paper.
Another natural direction for the research would be to investigate the generalization of different geometric properties in finite metric spaces. We only mention two of them here.
In elementary geometry, two distinct lines intersect in at most one point. If the same holds for a betweenness structure, then it is called \emph{geometric}.
Further, a betweenness structure is said to be \textit{Euclidean} if it is induced by a finite set of points of the $d$-dimensional Euclidean space for some nonnegative integer $d$. We note that every Euclidean betweenness structure is embeddable into the Euclidean plane.

It is easy to see that all Euclidean betweenness structures are regular, but not all regular betweenness structures are Euclidean: the simplest counterexample would be the betweenness structure induced by the star with $3$ leaves.
Similarly, all Euclidean betweenness structures are geometric but the reverse is false, as demonstrated by $\mathcal{C}_4$.
A more complicated example shows that even regularity and geometricity combined are not enough to guarantee the Euclidean property. Hence, characterization of Euclidean betweenness structures in a purely combinatorial way remains an interesting open problem.

\section*{Acknowledgment}
We are very grateful to Andreas Dress for helping us navigate through the literature on phylogenetic trees and metric decomposition theory of finite metric spaces.

\bibliographystyle{elsarticle-num}
\bibliography{fms_bib}

\end{document}